  \RequirePackage{lineno} 
\documentclass[12pt]{article}
\usepackage{amssymb}

\usepackage{amssymb,amsfonts}

\usepackage{pdflscape}
\usepackage{graphicx}
\usepackage{mathtools}
\usepackage{amsthm}
\usepackage{setspace}
\usepackage{diagbox}
\usepackage{adjustbox}
\usepackage{array}
\usepackage{multirow}
\usepackage{booktabs}
\usepackage[all,arc]{xy}
\usepackage{enumerate}
\usepackage{mathrsfs}
\usepackage{setspace}
\newtheorem{lemma}{Lemma}
\newtheorem{theorem}{Theorem}
\newtheorem{definition}{Definition}

\theoremstyle{definition}

\theoremstyle{remark}

\makeatletter
\let\c@equation\c@thm
\makeatother
\numberwithin{equation}{section}

\title{Formation of coalition structures as a  non-cooperative game}

\author{Dmitry Levando \thanks{ This paper comes as   a further development of the 3-rd and the 4-th chapters  of my PhD Thesis ``Essays on Trade and Cooperation'' at Ca Foscari University, Venezia, Italy. Acknowledge for comments and suggestions at different stages of the project:  Nick Baigent, Phillip Bich, Jean-Marc Bonnisseau, Alex Boulatov, Emiliano Catonini, Giulio Codognato, Sergio Currarini, Luca Gelsomini,  Izhak Gilboa, Olga Gorelkina, Piero Gottardi, Roman Gorpenko, Eran Hanany, Natalya Kabakova, Anton Komarov, Ludovic Julien, Alex Kokorev, Alex Larionov, Dmitry Makarov, Francois  Maniquet, Igal Miltaich, Stephane Menozzi,   Roger Myerson, Olga Pushkareva, Ariel Rubinstein, Marina Sandomirskaya,  Konstantin Sonin, Daria Schepetova,  William Thompson, Simone Tonin, Dimitrios Tsomocos, Eyal Winter, Shmuel Zamir. Special thanks to Nadezhda Likhacheva, Fuad Aleskerov,  Lev Gelman, Mark Kelbert, O.P., Shlomo Weber. Many thanks for  advices and for discussion to participants of  SI\&GE-2015, 2016 (an earlier version of the title was "A generalized Nash equilibrium"), CORE 50  Conference, CEPET 2016  Workshop, Games 2016 Congress, Games and Applications 2016 at Lisbon, Games and Optimization at St Etienne, Combinatorial Games Colloquium at Lisbon (2017). Previous version of the paper is available at www.arxiv.org.  All possible mistakes are only mine. This research did not receive any specific grant from funding agencies in the public, commercial, or not-for-profit sectors. \newline E-mail for correspondence: dlevando (at)  hse dot ru, d\_levando (at) hotmail dot com}}
\date{} 

\begin{document}

\maketitle
\begin{abstract}

Traditionally  social sciences are interested in structuring people in multiple groups based on their individual preferences. This paper suggests an approach to this problem in the framework of a non-cooperative game theory.


 Definition of a suggested finite game includes  a family of nested simultaneous non-cooperative  finite games   with intra- and inter-coalition externalities.  In this family,  games differ  by  the  size of   maximum coalition, partitions and  by coalition structure formation rules. 
 
A result of  every game  consists of  partition of   players into coalitions and a payoff profile for every player.  
 Every  game in the family has an equilibrium in mixed strategies with possibly  more than one coalition.  
 The results of the game  differ   from  those conventionally discussed in cooperative game theory, e.g. the Shapley value,  strong Nash, coalition-proof equilibrium,  core, kernel, nucleolus. 
 
 We discuss  the following applications of the new game:  cooperation as an allocation in one coalition,  Bayesian games, stochastic games and construction of a non-cooperative criterion  of coalition structure stability for studying  focal points.  
 \end{abstract}



 



\section{Introduction}


There is a conventional dichotomy  in game theory applications,  the cooperative game theory (CGT) versus  the non-cooperative game theory (NGT). CGT deals with coalitions as elementary items, NGT deals with strategic individual behavior.

However, the picture is not so smooth. CGT traditionally disregards issues of    strategic  interactions  between individual players,    with players in the same coalition and  players   in other coalitions.    CGT substitutes an individual gain by a value of a coalition,  a total and not personified gain of everybody.  As a result,  value  of a coalition does not depend on an allocation of all other players, and  construction of an individual gain   requires additional arrangements. 
 From another side NGT with an individual impact on  an equilibrium overlooks   structural aspects of player{'}s partition into groups or into coalition structures.\footnote{Coalition structure terminology was used by Aumann and Dreze (1976).} 
 
 The two theories even pursue different goals,  their outcomes vary in many ways, and are hard to compile.    CGT concentrates on socially desirable allocations of players, while efficiency is defined in aggregate terms, rather than  individual, as in the standard definition of Pareto-efficiency. 
  At the same time NGT aspires towards equilibrium conditions, which later could  be Pareto-improved. 
  
   Thus it is not surprising that  there are questions  beyond the range of concern for any of the two theories.
 The present study   concentrates on   games, where players can be allocated to different coalitions in an equilibrium,  still preserving traditional for NGT   individual strategies and individual payoffs. This results in two types of externalities from every player,  intra-coalition and inter-coalition ones.
 
Consider   two  similar examples:  a voluntarily division of   participants into paintball teams or a voluntarily division of a class into studying groups. 
Before they are engaged into a team/a group every  player/student makes a decision from self-interest considerations, concerning with whom he/she is,  and how all others are allocated. 
A decision  includes a choice of a  preferable  coalition structure, and what to do after it is  formed. 

  In these  examples the participants know about  commonly accepted mechanisms of coalition formation before any activity starts. These mechanisms resolve possible conflicts between individual choices. 
  Members of the team/the group may  perform  different actions  inside a team/a group\footnote{Indefinite article is used as a formed  team/group may not coincide with   preferable.  See an example in Section 2.} formed.\footnote{This is different from an approach within cooperative game theory, where players from a coalition  perform the same action.}

  Intra- and inter-coalition externalities for the  paintball game are clear. 
  Let us discuss two type of externalities for studying group game. There are  maybe  intra-coalition externalities in studying groups in a form of mutual assistance.  Inter-coalition externalities for studying groups could be  a mutual noise, or a slow WiFi connection from an external service.
 
The  common features of  the examples include multi-coalition frameworks, formed from self-concerned (self-interested) behavior of everybody, and  two types of externalities. Thus both examples occur in an area between the existing  cooperative and non-cooperative game  theories.

 
In this paper  `coalition structure', or  `partition' \footnote{Existing literature uses both terms.}  for short, is a disjoint union of non-overlapping subsets from a set of players. A group, or a coalition, is an element of a coalition structure or  of a  partition.\footnote{the same.}



 The examples above may seem to be related to cooperation, as suggested by Nash (1953). 
 Now  his suggestion  is known   as the   Nash Program (Serrano, 2004), and   cooperation   is understood as an activity inside  a group with positive externalities between players. However both examples above are more complicated than the existing theories suggest. For example, due to a multi-coalition framework and two types of externalities. 
 
The best analogy for the difference of  the Nash Program  from    the   research of this paper is the  difference between a partial and  a general strategic equilibrium analysis\footnote{meaning, general equilibrium framework with strategic market games of Shapley and Shubik.} in economics. The former isolates a  market  and ignores cross-market interactions, the latter explicitly studies cross market interactions from individual  strategic actions of self-interested traders.



 The  novelty of the paper consists in: a  construction of a \textit{family} of nested non-cooperative games such that every game has an equilibrium in mixed strategies. Based on the game, the paper presents a non-cooperative criterion for coalition structure stability, that leads to a non-cooperative  reconsideration of focal points of Schelling. These results  reveal new properties of well-known games, and   ask new research questions.\footnote{ Recently   Dasgupta and Maskin (2016) wrote that 
'game theory is the study of strategic interaction: how a person should and will behave if her actions affect others and their actions affect her`. The novelty of this paper is that  players can be allocated in different coalitions still making externalities for all others. }   
 Using   a game similar to Prisoner{'}s Dilemma the paper shows that  an efficient and not equilibrium outcome does not imply allocation of players in one coalition.  Further reading to the literature on this subject is available in the Discussion section of this paper.

 The paper addresses only construction and equilibrium issues of a non-cooperative game, but does not study  efficiency of an equilibrium concept.   


The  paper has the following structure: Section 2 presents an example based on Prisoner{'}s Dilemma and demonstrates that there is actually no explicit cooperation in it. Then  Section 3 follows  with the main model.   Section 4 contains a discussion of the result and a comparison with existing literature. Section 5 contains a formal definition of a sufficient criterion for cooperation. All  further sections deal with  reconsidering  well-known games,  Bayesian games (Section 6), stochastic games (Section 7),  and demonstrating their new properties, including a  construction of non-cooperative criterion of coalition structure stability and reconsidering focal points  of Shelling (Section 8).  The final section of this paper is a Conclusion. 

\section{Modified  Prisoner{'}s Dilemma  (PD)}
The example  demonstrates that there is no explicit cooperation in   the   standard PD game, if to think about a  cooperation  as an allocation of players in the same coalition. Formal analysis of cooperation is in
Section 5.\footnote{This section contains an example for two players, what may make to think that the players care about coalitions, not about coalition structures. The suggested mechanism is constructed to deal with  coalition structures. A game of two ``introvert'' players further in this paper demonstrates importance of coalition structures.  An example, where coalition structure matter explicitly is in Section 6, Corporate lunch game. }
 
 Consider a game of two players $i=1,2$, where each can choose   to be alone and has two activity levels:  L(ow)$_{i,alone} $  or  H(high)$_{i,alone} $, $i = 1,2$. The players can make only one coalition structure $\{\{ 1\} ,\{ 2\} \}$.  Payoffs for this case are in Table \ref{PDstandard}.
 Every cell contains a payoff profile for both players and an allocation of the players over coalitions.  The unique equilibrium $ {(-2;-2)^{*}}$, marked with one star $^{*}$, is inefficient. But  the desirable  and non-equilibrium payoff outcome $(0;0)$ is efficient. Hence the game reproduces payoff properties of the Prisoner{'}s Dilemma with an  explicit allocation of players over coalitions structures. 
 
 \begin{table}[htp]
\caption{Payoff for the standard Prisoners Dilemman}   
 \begin{center}
\begin{tabular}{|c|c|c|} \hline
 & $L_{2,alone}$ & $H_{2,alone} $  \\  \hline
$L_{1,alone}$ &
 \begin{tabular}{c} \underline{(0;0)}  \\ $ \{\{1\},\{ 2\} \}$ \end{tabular}
&  \begin{tabular}{c}  (-5;3) \\ $\{\{1\},\{ 2\} \}$  \end{tabular} \\ \hline
$H_{1,alone}$ 
&  \begin{tabular}{c}  (3;-5) \\ $ \{\{1\},\{ 2\} \}$ \end{tabular}
&  \begin{tabular}{c} $ {(-2;-2)^{*}}$ \\ $  \{\{1\},\{ 2\} \}$ \end{tabular}
\\ 
\hline 
  \end{tabular}
\end{center}
\label{PDstandard}
\end{table}%

Let the game be more complicated, and  the players can choose to be either alone or together, and  as above to choose between two   activity levels for every case.   Then set of strategies $S_i(K=2)$ for player  $i$ is:
 
 $$S_i(K=2)= \Big\{(L_{1,alone},H_{1,alone}), (L_{1,together},H_{1,together})  \Big\}.$$
 
  Parameter $K=2$ indicates that   maximum available coalition size for the players is  2. For simplicity we  assume  the simplest rule for coalition structure formation:  if at least one of the players  chooses to be alone, then a coalition of  two cannot be formed. 
 Payoff matrix for this game is in Table  \ref{PD2}.
 
 \begin{table}[htp]
\caption{Payoff for an extension of Prisoner{'}s Dilemma.  }
\begin{center}
\begin{tabular}{|c|c|c||c|c|} \hline
 & $L_{2,{alone}}$ & $H_{2,{alone}} $& $L_{2,{together}} $& $H_{2,{together}} $ \\  \hline
$L_{1,{alone}}$ &
 \begin{tabular}{c} \underline{(0;0)}  \\ $ \{\{1\},\{ 2\} \}$ \end{tabular}

&  \begin{tabular}{c}  (-5;3) \\ $\{\{1\},\{ 2\} \}$  \end{tabular}

& \begin{tabular}{c} \underline{(0;0)} \\  $  \{\{1\},\{ 2\} \}$ \end{tabular}
& \begin{tabular}{c}  (-5;3) \\ $  \{\{1\},\{ 2\} \}$\end{tabular}  \\ \hline
$H_{1,{alone}}$ 
&  \begin{tabular}{c}  (3;-5) \\ $ \{\{1\},\{ 2\} \}$ \end{tabular}
&  \begin{tabular}{c} $ {(-2;-2)^{*,**}}$ \\ $  \{\{1\},\{ 2\} \}$ \end{tabular}
&  \begin{tabular}{c}  (3;-5) \\ $  \{\{1\},\{ 2\} \}$ \end{tabular}
& \begin{tabular}{c}  ${ (-2;-2)}^{**}$ \\ $  \{\{1\},\{ 2\} \}$\end{tabular} \\ 
\hline \hline
$L_{1,{together}}$ 
&  \begin{tabular}{c}  \underline{(0;0)} \\ $  \{\{1\},\{ 2\} \}$\end{tabular}
&  \begin{tabular}{c}  (-5;3) \\ $  \{\{1\},\{ 2\} \}$\end{tabular}
&  \begin{tabular}{c}  $  \underline{(0;0)} $ \\ $  \{1, 2 \}$\end{tabular}
&  \begin{tabular}{c}  $(-5;3)$  \\ $  \{1,2\}$  \end{tabular}\\ \hline
$H_{1,{together}}$ 
&  \begin{tabular}{c}  (3;-5) \\ $  \{\{1\},\{ 2\} \}$ \end{tabular}
&  \begin{tabular}{c}  ${ (-2;-2)}^{**}$ \\ $  \{\{1\},\{ 2\} \}$ \end{tabular}
& \begin{tabular}{c}   $(3;-5) $  \\ $  \{1,2\}$ \end{tabular}
&  \begin{tabular}{c}  ${ (-2;-2)}^{**}$ \\ $  \{1, 2\}$   \end{tabular}\\ \hline  
 \end{tabular}
\end{center}
\label{PD2}
\end{table}%

 Every cell in Table \ref{PD2} contains a payoff profile and a coalition structure. The upper-left part of  Table \ref{PD2} corresponds to strategies $L_{i,alone} $ and  $H_{i,alone} $, $i=1,2$, and coincides with  the previous game in Table \ref{PDstandard}.
 The previous equilibrium is marked by  one star $^{*}$, the newly discovered  equilibria are marked by  two stars $^{**}$. It is clear that all new equilibria are still inefficient, but   belong to  different coalition structures.  An increase in a number of possible coalition structures enriches a set of equilibria of the game and   a set of equilibrium coalition structures, but not efficiency.

  In the coalition structure $\{ \{1 \} , \{2 \}\}$ two  players have  externalities within coalitions, which will be further referred to as inter-coalition externalities. But in the coalition structure $\{ 1 , 2 \}$ the same players  experience  between externalities,  will be further referred to  as intra-coalition externalities. 
  
Formation of   $\{1,2 \}$ requires that two  players choose to be together, whatever they plan to do inside the coalition structure.  It is  impossible to implement  $\{1,2 \}$   only from a choice of   one player,   given  coalition structure formation rule. Inside this coalition structure both players can deviate simultaneously from the efficient outcome to the equilibrium. 

The  feature  of the example is that there is no instruction for players on how to deviate:  unilaterally, independently, or  within only one coalition, as in  the cooperative game theory.  \textit{The sets of individual strategies anticipate all combinations of possible deviations for both players.} A formal model in Section 2 generalizes this idea.
  
 A resulting coalition structure is formed  as  \textit{a result of  choices from all players} for a given coalition structure formation rule. Exogenous rules resolve conflicts  between choices of players becoming  components of a game.

  Table \ref{PD2} contains four efficient outcomes: $(0;0)$ located in different   coalition structures: this means that   by observing only the payoff profiles $(0;0)$ we can{'}t conclude   whether   players  are in one coalition or   not. This differs from a common   understanding that the Prisoner{'}s Dilemma  can serve as  an example of   unsuccessful  cooperation. If to consider a  cooperation as  an allocation of players in one coalition,   we observe, that information about payoff profile is not enough to assert cooperation existence.

 An   allocation of players in one coalition can be fixed by addressing to the game of  ``extrovert'' players.  Two players have preferences over coalition structures and prefer to be together. Let $\epsilon >0$ be an additional individual  gain   (a corporate gain) by contrast from being separate, but only if the grand coalition $\{1,2 \}$ is formed. Table  \ref{PD3} contains payoffs for this game.
  {\tiny{\begin{table}[t]
\caption{Payoff for two extrovert players who  prefer to be together. Uniqueness of  an equilibrium is fixed.}
\begin{center}
\begin{tabular}{|c|c|c||c|c|} \hline 
 & $L_{2,alone}$ & $H_{2,alone} $& $L_{2,together} $& $H_{2,together} $ \\  \hline
$L_{1,alone}$ &
 \begin{tabular}{c} \underline{(0;0)}  \\ $ \{\{1\},\{ 2\} \}$ \end{tabular}

&  \begin{tabular}{c}  (-5;3) \\ $\{\{1\},\{ 2\} \}$  \end{tabular}

& \begin{tabular}{c} \underline{(0;0)} \\  $  \{\{1\},\{ 2\} \}$ \end{tabular}
& \begin{tabular}{c}  (-5;3) \\ $  \{\{1\},\{ 2\} \}$\end{tabular}  \\ \hline
$H_{1,alone}$ 
&  \begin{tabular}{c}  (3;-5) \\ $ \{\{1\},\{ 2\} \}$ \end{tabular}
&  \begin{tabular}{c}  ${(-2;-2)}^{*}$ \\ $  \{\{1\},\{ 2\} \}$ \end{tabular}
&  \begin{tabular}{c}  (3;-5) \\ $  \{\{1\},\{ 2\} \}$ \end{tabular}
& \begin{tabular}{c}  ${ (-2;-2)} $\\ $  \{\{1\},\{ 2\} \}$\end{tabular} \\ 
 \hline \hline
$L_{1,together}$ & 
\begin{tabular}{c} \underline{(0;0)} \\ $ \{\{1\},\{ 2\} \}$\end{tabular}
& \begin{tabular}{c} (-5;3) \\ $  \{\{1\},\{ 2\} \}$ \end{tabular}
& \begin{tabular}{c} $  \underline{(0+\epsilon;0+\epsilon)} $ \\  $\{1, 2 \}$ \end{tabular}& \begin{tabular}{c} $(-5+\epsilon;3+\epsilon)$ \\  $\{1,2\}$ \end{tabular} \\ \hline
$H_{1,together}$ 
& \begin{tabular}{c} (3;-5) \\ $ \{\{1\},\{ 2\} \}$\end{tabular} 
&  \begin{tabular}{c} (-2;-2)  \\ $  \{\{1\},\{ 2\} \}$ \end{tabular} & 
\begin{tabular}{c}
$(3+\epsilon;-5+\epsilon) $ \\ $ \{1,2\}$ \end{tabular} &  \begin{tabular}{c}  ${ (-2+\epsilon;-2+\epsilon)}^{**}$ \\ $ \{1,2 \}$  \end{tabular} \\ \hline  
 \end{tabular}
\end{center}
\label{PD3}
\end{table}%
}}
Now the equilibrium marked with $^{**}$ is unique  and  appears  when  both players   
choose to be together. However it   is non-efficient again.  

If both players are  ``introverts'', then each prefers to be alone. If a desirable coalition structure, $\{ \{1\},\{ 2\}\}$, is formed, then every player obtains an additional markup  $\delta >0$.    Payoff matrix for this game is in Table \ref{introverts}. The additional option of being  together does not change  the equilibrium in comparison to the aforementioned game, where strategies are $(L_{i,alone},H_{i,alone})$, $i=1,2$.

    {\tiny{\begin{table}[h]
\caption{Payoff for two ``introvert'' players, who prefer to be alone.  Uniqueness of  an equilibrium is fixed.}
\begin{center}
\begin{tabular}{|c|c|c||c|c|} \hline 
 & $L_{2,alone}$ & $H_{2,alone} $& $L_{2,together} $& $H_{2,together} $ \\  \hline
$L_{1,alone}$ &
 \begin{tabular}{c}  $\underline{(0+\delta;0+\delta)} $ \\ $ \{\{1\},\{ 2\} \}$ \end{tabular}

&  \begin{tabular}{c}  $(-5+\delta;3+\delta)$ \\ $\{\{1\},\{ 2\} \}$  \end{tabular}

& \begin{tabular}{c} $\underline{(0;0)} $\\  $  \{\{1\},\{ 2\} \}$ \end{tabular}
& \begin{tabular}{c}  $(-5;3)$ \\ $  \{\{1\},\{ 2\} \}$\end{tabular}  \\ \hline
$H_{1,alone}$ 
&  \begin{tabular}{c} $ (3+\delta;-5+\delta)$ \\ $ \{\{1\},\{ 2\} \}$ \end{tabular}
&  \begin{tabular}{c}  ${(-2+\delta;-2+\delta)}^{*,**}$ \\ $  \{\{1\},\{ 2\} \}$ \end{tabular}
&  \begin{tabular}{c}  (3;-5) \\ $  \{\{1\},\{ 2\} \}$ \end{tabular}
& \begin{tabular}{c}  ${ (-2;-2)} $\\ $  \{\{1\},\{ 2\} \}$\end{tabular} \\ 
 \hline \hline
$L_{1,together}$ & 
\begin{tabular}{c} \underline{(0;0)} \\ $ \{\{1\},\{ 2\} \}$\end{tabular}
& \begin{tabular}{c} (-5;3) \\ $  \{\{1\},\{ 2\} \}$ \end{tabular}
& \begin{tabular}{c} $  \underline{(0;0)} $ \\  $\{1, 2 \}$ \end{tabular}& \begin{tabular}{c} $(-5;3)$ \\  $\{1,2\}$ \end{tabular} \\ \hline
$H_{1,together}$ 
& \begin{tabular}{c} (3;-5) \\ $ \{\{1\},\{ 2\} \}$\end{tabular} 
&  \begin{tabular}{c} (-2;-2)  \\ $  \{\{1\},\{ 2\} \}$ \end{tabular} & 
\begin{tabular}{c}
$(3;-5) $ \\ $ \{1,2\}$ \end{tabular} &  \begin{tabular}{c}  ${ (-2;-2)}$ \\ $ \{1,2 \}$  \end{tabular} \\ \hline  
 \end{tabular}
\end{center}
\label{introverts}
\end{table}%
}}

Consider a case, where the  players are different: player $i=1$ is ``extrovert'', but player $i=2$ is ``introvert''. Table \ref{PD4}  is the payoff matrix for this game. Following the rule of unanimous agreement to form a coalition the grand coalition, $\{ 1,2\}$, can not be formed, and equilibrium strategy profile for $i=2$ does not change in comparison to the first example. 
For player $i=1$  an equilibrium strategy profile is a mixed strategy with equal weights over two pure strategies $H_{1,alone}$ and $H_{1,together}$.   Final coalition structure  $\{\{ 1\} , \{ 2\} \}$ is not affected by this randomization  and another  coalition structure  is unfeasible. Players have inter-coalition externalities and inefficient outcome. Hence   the  mixed strategies induce
inter-coalition externalities and two coalitions in an equilibrium.  Games with such properties  are not described in existing literature.

{\tiny{\begin{table}[t]
\caption{Payoff for a game when player 1 is "extrovert" and player 2 is "introvert".  }
\begin{center}
\begin{tabular}{|c|c|c||c|c|} \hline 
 & $L_{2,alone}$ & $H_{2,alone} $& $L_{2,together} $& $H_{2,together} $ \\  \hline
$L_{1,alone}$ &
 \begin{tabular}{c} \underline{$(0;0+\delta)$}  \\ $ \{\{1\},\{ 2\} \}$ \end{tabular}

&  \begin{tabular}{c}  $(-5;3+\delta)$ \\ $\{\{1\},\{ 2\} \}$  \end{tabular}

& \begin{tabular}{c} \underline{$(0;0+\delta)$} \\  $  \{\{1\},\{ 2\} \}$ \end{tabular}
& \begin{tabular}{c}  $(-5;3+\delta)$ \\ $  \{\{1\},\{ 2\} \}$\end{tabular}  \\ \hline
$H_{1,alone}$ 
&  \begin{tabular}{c}  $(3;-5+\delta)$ \\ $ \{\{1\},\{ 2\} \}$ \end{tabular}
&  \begin{tabular}{c}  ${(-2;-2+\delta)}^{*,**}$ \\ $  \{\{1\},\{ 2\} \}$ \end{tabular}
&  \begin{tabular}{c}  $(3;-5+\delta)$ \\ $  \{\{1\},\{ 2\} \}$ \end{tabular}
& \begin{tabular}{c}  ${ (-2;-2+\delta)} $\\ $  \{\{1\},\{ 2\} \}$\end{tabular} \\ 
 \hline \hline
$L_{1,together}$ & 
\begin{tabular}{c} \underline{$(0;0-5+\delta)$} \\ $ \{\{1\},\{ 2\} \}$\end{tabular}
& \begin{tabular}{c} $(-5;3+\delta)$ \\ $  \{\{1\},\{ 2\} \}$ \end{tabular}
& \begin{tabular}{c} $  \underline{(0+\epsilon;0)} $ \\  $\{1, 2 \}$ \end{tabular}& \begin{tabular}{c} $(-5+\epsilon;)$ \\  $\{1,2\}$ \end{tabular} \\ \hline
$H_{1,together}$ 
& \begin{tabular}{c} $(3;-5+\delta)$ \\ $ \{\{1\},\{ 2\} \}$\end{tabular} 
&  \begin{tabular}{c} $ (-2;-2+\delta)^{**} $ \\ $  \{\{1\},\{ 2\} \}$ \end{tabular} & 
\begin{tabular}{c}
$(3+\epsilon;-5) $ \\ $ \{1,2\}$ \end{tabular} &  \begin{tabular}{c}  ${ (-2+\epsilon;-2)} $ \\ $ \{1,2 \}$  \end{tabular} \\ \hline  
 \end{tabular}
\end{center}
\label{PD4}
\end{table}%
}}

   
   We can summarize the section: expansion of  traditional non-cooperative game theory   for  coalition structures formation  requires additional arrangements.  This includes: 
   \begin{enumerate}
\item  
to be aware of a    maximum available coalition size;   
\item to enumerate possible coalition structures for it;  
\item individual strategy sets should be related to possible coalition structures;  
\item   to embed  some mechanism to resolve conflicts between players.  
\end{enumerate}
The next section  presents a formal model for   a non-cooperative game with coalition structure formation.
      
\section{Formal setup  of the  model}
This section presents a construction of the main model. 
Nash (1950, 1951) suggested a non-cooperative game, that  consists of a set of players $N$, with a general element $i$, sets of  individual finite strategies $S_i,  i \in N$, and individual payoffs, $U_i(s)$, defined for every strategy profile $s \in S=\times_{i \in N}S_i$,  as a mapping from a set of all strategies $S$ into a bounded set $U_i  \colon S \mapsto \mathbb{R}$, $U_i(s) < \infty $. Thus every  strategy profile  $s  \in S$ is assigned a vector of payoffs $\Big(U_i(s) \Big)_{i \in N} \in \mathbb{R}^N_+$, where $u_i(s)$ is an individual payoff of $i$ for a strategy profile $s$.  However, the non-cooperative game introduced by Nash says nothing about coalition structure formation.

The   game suggested in this paper    follows   similar logic of Nash, but with a preliminary modification  for    a set of   strategy profiles. Every player has a set of strategies for every relevant coalition structure. A set of all strategy profiles  $S$, constructed as a Cartesian (direct) product of all individual strategies,  is  portioned  into non-overlapping  domains, as in Section 2. Every  constructed domain is a set of strategy profiles. Every strategy profile inside a  domain   is assigned  the same coalition structure and a profile of individual payoffs.  
This approach reproduces the method  of  state-contingent payoffs for Arrow-Debreu securities.   

  In this section we introduce   components for a game,  how to assemble them into a game  and now to construct  an equilibrium concept for the game. The common feature of the components is their nested structure, that enables us  to construct a nested family of games.
  
\subsection{ Players and partitions} 

Consider  a set of agents $N$, with a general element $i$, and denote a size of  $N$  as  $\#N$, a  finite integer, $2 \le \#N < \infty$. If there is no ambiguity we will use $N$ also for a size of the set of players.

Introduce   a parameter $K$, taking values from 1 to $N$. It  has two equivalent interpretations:   maximum coalition size and   maximum number of deviators. 
 From one side,  in any available coalition structure  there is no coalition with maximum size greater than  $K$. From another,   at most   $K$  agents are required to dissolve any coalition  in any such coalition structure. 
 
Let $n(k)$   be a subset of players with a size no more than $K$, $n(k) \subset N$, $\#n(k) \le K$. \textit{A subset $n(k)$ is not a coalition, but a subset of players from $N$ with an upper restriction on a number of elements.} After a coalition structure is formed the players from  $n(k)$ may appear in different  final coalitions. 

Every  fixed  $K$  induces a family of coalition structures  (or a family of partitions)   $\mathcal{P}(K)$: 
\begin{multline}
\mathcal{P}(K) =  \Big\{P \colon P=\{g_j \colon g_j \subset N, \ \  \#g \le K , \ \ \sqcup_j g_j = N  \}\Big\}  \end{multline}

An element $P=\{g_j \}$ of  $\mathcal{P}(K)$ is a partition of  $N$ into coalitions, where $j$ is an index of  a coalition in $P$.
Every  coalition has an upper boundary on its  size, $\#g_j \le K$.  The notation  $\sqcup_j g_j = N$ means that a player may participate only in one coalition: $    \forall j_1 \neq j_2,   \Rightarrow    {g}_{j_1} \cap {g}_{j_2} = \emptyset  $. Size of $\mathcal{P}$ is described by a partition function $p(N,K)$ ( reference). 

If we increase $K$ by one, then we need    to add  new coalition structures from  $\mathcal{P}(K+1) \setminus \mathcal{P}(K) $, and finally we obtain the nested families of partitions: 
  $$\mathcal{P}(K=1) \subset \ldots \subset \mathcal{P}(K) \subset  \ldots  \subset \mathcal{P}(K=N).$$  The bigger   $K$ is , the more  coalition structures (or partitions) are  included into consideration.  The grand coalition, i.e. a coalition of size $N$, belongs    to the family $\mathcal{P}(K=N) $, only.
  
\subsection{  Strategies }
Take a fixed $K$ and consider  a  relevant family of coalition structures $\mathcal{P}(K)$.    
     For every partition $P$ from $ \mathcal{P}(K)$ agent  $i$ has a finite strategy set $S_i (P)$ and a set of strategies of $i$ for $\mathcal{P}(K)$ is $S_i (K)=\cup_{ P \in \mathcal{P}(K)} S_i(P)$.\footnote{Finite strategies are  used as  in  Nash (1950).}  
\begin{definition}[individual strategy set]
   A set of strategies of agent $i$ for a family of coalition structures    $\mathcal{P}(K)$ is
$$
S_i(K) = \Big\{s_i   \colon s_i \in \cup_{ P \in \mathcal{P}(K)} S_i(P) \Big\}
$$ with a general element $s_i(K)$ or $s_i$ for simplicity. 
\end{definition}
Agent $i$ chooses $s_i$ from $S_i(K)$.  All agents make their choices simultaneously. \textit{ A choice  is a choice of a \textbf{desirable} partition and an action for this partition}.\footnote{A desirable partition may not realize due to  conflicts  in individual choices. It was demonstrated with examples in the previous section. A sequential game for dynamic formation coalition structures is in progress with Marc Kelbert and Olga Pushkareva. Both versions of the game   need a mechanism for conflict resolution.} If we increase    $K$ by one, then we need to construct additional strategies only for the newly available coalition structures from  $\mathcal{P}(K+1) \setminus \mathcal{P}(K) $. This makes strategy sets for different $K$'s to   be nested: 
  $$S_i(K=1) \subset \ldots \subset S_i(K) \subset  \ldots  \subset S_i(K=N).$$ 
   
 
We construct   sets of strategies  of all players   from individual strategies in the usual way, $S(K)=\times_{i \in N} S_i (K)$, and it  inherits  the nesting property: 
   $$S(K=1) \subset \ldots  \subset S(K=N).$$

A choice of all players, a strategy profile, is $s(K)=\Big(s_1(K), \ldots,s_N(K)\Big)\in S(K)$, or if there is no ambiguity   $s   \equiv s(K)$. 

 One may ask a question,  is it possible  to make a choice in  two stages:  first,  all players choose a partition, second, everybody chooses a coalition.
Such reformulation of the game leads to   the  generalization of   strategic equilibrium, introduced by Mertens (1995), also in Hillas,  Kohlberg (2002),  that goes beyond the scope of this paper. 


\subsection{ Coalition structure mechanism}
  
  For every  $K$  we have a family of coalition structures and a set of strategy sets. However, individual choices of agents may be inconsistent. To overcome this deficiency  we define a coalition structure formation mechanism  $\mathcal{R}(K)$, (a mechanism  or a rule for short) to resolve possible conflicts.\footnote{They can be social norms or social institutions, which people take for granted.}  
 The rule portions   $S(K)$ into  non-overlapping coalition structure specific domains 
 $S(K)=\sqcup_{P \in \mathcal{P}(K)} S (P)$,  $S (P) = \times_{i \in N} S_i(P)$, where disjoint union $\sqcup$ means that domains for different coalition structures do not overlap. An example was in Section 2. 
     
   A unique  coalition structure $P$ from $ \mathcal{P}(K)$    is assigned for  every  strategy profile $s$ from $S(K)$.  A set of such strategy profiles makes a coalition structure  strategy domain $S(P)=\{s=(s_1,\ldots,s_i,\ldots,s_N) \}$, which  may not be     convex,  (see examples in Section 2), $s_i \in S_i(K)$.  It is clear that $S(K)=\sqcup_{ P \in \mathcal{P}(K) }S(P)$. For some players a realized coalition structure $P$ may not coincide with a chosen one.

\begin{definition}[a coalition structure formation mechanism] For every $K$  \textbf{a coalition structure formation mechanism}  $\mathcal{R}(K)$  is a  set of measurable mappings such that:
 \begin{enumerate}
\item A domain of $\mathcal{R}(K)$  is a set of all strategy profiles  of  $S(K)$. 
\item A range of $\mathcal{R}(K)$   is a finite  number of subsets  $S(P) \subset S(K)$, $P \in \mathcal{P}(K)$. Every $S(P)$ is a strategy set for only one coalition structure $P$.
\item   $\mathcal{R}(K)$ divides $S(K)$  in such a way that the union of all $S(P)$ constructs the original set $ S(K) = \sqcup_{ P \in \mathcal{P}(K)}S(P)$. 
\end{enumerate}
Formally the same:
 \begin{multline*}\mathcal{R}(K) \colon 
S(K)=\times_{i \in N}S_i(K)  \mapsto    \Big(S(P), P \in \mathcal{P}(K) \Big) \colon \sqcup_{P \in \mathcal{P} (K)} S(P)
 \end{multline*}   
\end{definition}

There are two ways to construct $S(K)$: in terms of initial individual strategies, $S(K)=\times_{i \in N}S_i (K)$, or  in terms of   realized partition   strategies,  $S(K) =\sqcup_{ P \in \mathcal{P}(K)} S(P)$. 

Holt and Roth (2004) suggested to   include  a mechanism  to eliminate a distinction between non-cooperative and cooperative games: `One trend in modern game theory, often referred as \textit{Nash program}, is to erase this distinction by including any relevant enforcement mechanism  in the model, so that all games can be modeled as non-cooperative'.  

  Informally, every coalition structure   becomes  itself  a   non-cooperative game: for  every coalition structure $P$ there is  a set of players, a  non-trivial set of strategies and  partition-specific payoffs, described  further. 
 
If $K$ increases we need to add   a mechanism for strategy sets  from   $S(K+1) \setminus S(K)$ only.   This supports  consistency of  coalition structure formation mechanisms for different $K$, and the  family of mechanisms is nested: 
 $$\mathcal{R}(K=1)  \subset \ldots \subset \mathcal{R}(K) \subset  \ldots  \subset \mathcal{R}(K=N).$$  
 
 It is important to note that  a role of a  mechanism is to   identify domains for coalition structures, but  not to deal with    efficiency  of results (see inefficient equilibria in Section 2).\footnote{Efficiency in non-cooperative and cooperative game theories is understood differently.} Studying efficiency requires studying  properties of all possible measurable transformations $\mathcal{R}(K)$, that goes far beyond this paper.  Discussion of  some properties of $\mathcal{R}(K)$ is in the next section. 
 
 \subsection{  Payoffs} 
 For every  $K$  we have a family of coalition structures, a set of strategy sets and coalition structure formation rules. And now  within every given coalition structure for every strategy profile we define a   payoffs profile.   Payoffs are considered as von  Neumann-Morgenstern utilities. 
 
  
  \begin{definition}[coalition structure specific individual payoff]For every  coalition  structure  $P$ from $\mathcal{P}(K)$   player $i$, $i \in N$, has a payoff function $$U_i (P) \colon S(P) \rightarrow \mathbb{R}_+,$$ such that    $U_i (P) $ is  bounded, $U_i    < \infty$.\end{definition}
 Definition of payoffs for all  feasible coalition structures generalizes the   definition above. 
   \begin{definition}[family of   coalition structure specific individual payoffs]
   For  given maximum coalition size $K$, family of coalition structures $\mathcal{P}(K)$,  set of strategy profiles  $S(K)$,   payoffs  of  player $i$ from $N$, make  a family:
  $$\mathcal{U}_i(K) = \Big\{U_i\Big(P \Big) \colon   P \in \mathcal{P}(K), S(K)=\sqcup_{P \in \mathcal{P}(K)} S(P)\Big\},
   $$ where $U_i\Big( P\Big)$ is a set of coalition structure specific payoffs for $i$. 
   \end{definition}
   Every coalition structure   has a
 strategy domain  $S(P)$,  and   a  set of payoffs, $\Big(U_i\Big(S(P)\Big)\Big)_{i \in N}$, defined on this domain.  Thus, every coalition structure   is a non-cooperative game itself.  
  
 An increment in $K$   increases a number of possible partitions and a set of  feasible strategies for every player. Thus we need  to add only payoffs for newly added coalition structures from   $\mathcal{P}(K+1) \setminus \mathcal{P}(K)$.   
Not surprisingly we  obtain a nested family of payoff functions:
 $$ \mathcal{U}_i(K=1) \subset  \ldots \subset  \mathcal{U}_i(K) \subset  \ldots \subset \mathcal{U}_i(K=N).$$ 

It is clear that in every coalition structure players can   have  both intra and inter coalition externalities. 

\subsection{ Game} 
For every K we have a family of coalition structures, strategy sets, coalition structure formation rules and families of payoffs and  now, we define a game.
\begin{definition}[ {a simultaneous coalition structure formation game}]
A non-cooperative game   for coalition structure formation with  maximum coalition size $K$ is  
$$\Gamma(K)=\Big\langle N, \Big\{K, \mathcal{P}(K),\mathcal{R}(K) \Big\}, \Big(S_i(K),\mathcal{U}_i (K)\Big)_{i \in N} \Big\rangle, $$  where
$\Big\{K, \mathcal{P}(K),\mathcal{R}(K) \Big\}$  - coalition structure formation mechanism (a social norm, a social institute),
$\Big(S_i(K),\mathcal{U}_i (K)\Big)_{i \in N}$ - properties of  players in $N$, (individual strategies and payoffs), such that: 
 $$\times_{i \in N} S_i(K)\buildrel \mathcal{R}(K) \over 
 \rightarrow \Big\{ S(P) \colon P \in \mathcal{P}(K) \Big\}
 \rightarrow  \Big\{\Big(\mathcal{U}_i(K)  \Big)_{i \in N} \Big\}_{ P \in \mathcal{P}(K)}. $$
\end{definition} 

If we omit the mechanism part of the game, i.e. take $K=1$,  and forget about the coalition structures, then  we obtain the traditional non-cooperative game of Nash. 
 It is clear that  for different $K$'s the constructed games   are nested. 
 
\begin{definition}[\textbf{family of games}]
A family of games is \textbf{ nested}  if :
$$
\mathcal{G}=\Gamma(K=1) \subset \ldots \subset \Gamma(K) \subset \ldots \subset \Gamma(K=N).
$$
 \end{definition}
 Nested games are important to study changes in an equilibrium with an increase in a number of deviators  or maximum coalition{'}s size $K$.
 
   
   \subsection{  Equilibrium} 
   Let $S_i (K)$ be a finite set of strategies of player $i$.
Let $\Delta_i(K)$ be a set of all mixed strategies  (probability measures or probability distributions) of player $i$, i.e. a space of probability measures over $S_i(K)$,  $$\Delta_i (K) = \Big\{\sigma_i (K) \colon \sum_{P \in \mathcal{P}(K)} \sum_{s_i \in S_i(K)} \sigma_i(s_i)=1 \Big\}.$$  Let   $\sigma_i (K)$ be a mixed strategy of $i$, $\sigma_i (K) \in \Delta_i (K)$.


Sets of mixed strategies for  players different from $i$ are defined in the standard way: $$\Delta_{-i} (K) = \Big\{\sigma_{-i}=\prod_{j \neq i} \sigma_j(K)  \colon  \forall j \neq i  \Big\}.$$
  
  Set of all strategies $S(K)$ has two representations, $S(K)=\times_{i \in N} S_i (K)$ and $S(K)=\sqcup_{P \in \mathcal{P}(K)} S(P)$. Thus  there are two parallel and equivalent ways to  construct expected utility.  
  
Expected utility of $i$ in terms of individual strategies $(S_i(K))_{i \in N}$ is 
$$
EU^{\Gamma(K)}_i \Big(\sigma_i(K),\sigma_{-i}(K)\Big) = 
\sum_{P \in \mathcal{P}(K)}\sum_{s_i \in S_i(P)} {U}_i(s_i,s_{-i}) \sigma_i(K) \sigma_{-i}(K),$$ and is defined in terms of realized coalition structures from $\mathcal{P}(K)$.


 \begin{definition}[\textbf{an equilibrium in a game $\Gamma(K)$ }]
There is a non-cooperative game for coalition structure formation $$\Gamma(K)=\Big\langle N, \Big\{K, \mathcal{P}(K),\mathcal{R}(K) \Big\}, \Big(S_i(K),\mathcal{U}_i (K)\Big)_{i \in N} \Big\rangle. $$
A mixed strategies profile  $\sigma^*(K)=\Big(\sigma^*_i(K)\Big)_{i \in N}$ is  an equilibrium  for a game $\Gamma(K)$ if for every subset $n(k)$ from $N$, with a size $1\le k \le K$, and  for every  player   $i \in n(k)$ a deviation from an equilibrium,     $\sigma_i(K) \neq \sigma^*_i(K)$,  does not generate an individual gain: 
 $$
 EU^{\Gamma(K)}_i \Big(\sigma^*_i(K),\sigma^*_{-i}(K)\Big) \ge EU^{\Gamma(K)}_i \Big(\sigma_i(K),\sigma^*_{-i}(K)\Big).
  $$
\end{definition}

The difference with Nash equilibrium is that an equilibrium must be satisfied not for every player $i$ from $N$, but for every player $i$ from every subset
$n(k)$ from $N$, where a  size   of  $n(k)$   does not exceed  $K$. Informally this means that any $k \le K$ deviators do not have restriction to move either separately or together. Players from $n(k)$ can deviate in any combinations, what is anticipated in their strategy sets in  construction of the game, 

Possible deviations from an equilibrium are   constructed differently from those used in the strong Nash equilibrium, in the core or in the coalition-proof equilibrium. 
This property of the simultaneous deviation  already appeared in the Prisoner{'}s Dilemma example, where if players  have an option to be together, then they must be aware of a deviation of both in an equilibrium. 
  \begin{theorem}
A game $\Gamma(K)$ has an equilibrium in  mixed strategies. 
\end{theorem}
A set of finite strategies has a simplex as a set of mixed strategies. Simplex is a compact set. Set of mixed strategies  of all players is a combination of finite number of simplex, what is again a compact set.
Continuous  utility function   function  reaches a maximum value on this simplex. Thus there is always a Nash equilibrium for finite strategies in $\Gamma (K)$.
Nash theorem with finite strategies is the partial case  for $K=1$.  The result of the paper is not in a trivial mathematical theorem, but in   introduction of a new game  to solve earlier untreatable problems with any number of deviators structured by groups. This new game   follows the same rules, as the game  suggested by Nash.

The equilibrium result is different from results of the cooperative game theory, where    games with empty cores may exist. The  suggested non-cooperative  equilibrium does not demand super-additivity, transferable / non-transferable utilities, axioms on payoffs or weights.  Existence of an equilibrium does not imply that it is  efficient. 

 The advantage of the theorem is that it can be generalized for any  relevant $K$.
 New insight of the result is not in complexity of mathematical objects, but in constructing a  more general structure of a game, which allow to ask new questions and to reconsider old games.

 \begin{theorem}
The family of games $\mathcal{G}=\{ \Gamma(K), K=1,2, \dots,N\}$ has equilibrium   mixed strategies profiles for every $K$:  $$\sigma^*(\mathcal{G})=\Big(\sigma^*(K=1), \ldots, \sigma^*(K)_{i \in N},\ldots, \sigma^*(K=N) \Big),$$  where $\sigma^*(K)=\Big(\sigma^*_i(K)\Big)_{i \in N}
$ is an equilibrium for a game $\Gamma(K)$.
\end{theorem}

The proof is simple.  A players $i${'}s finite set of pure strategies has simplex as a set of mixed strategies. Then a set of all strategies of a game $\Gamma(K)$ is a  convex subset of $\mathcal{R}^{p(N,K)}$, where $p(N,K)$ is partition function. Then by Brouwer{'}s fixed point theorem a continuous mapping from  a convex compact set  in a finite dimensional Euclidean space  into a convex subset  of the same  finite dimensional Euclidean space has a fixed point. 
  
Thus every game for a  coalition structure formation has an equilibrium for any number of deviators. It is clear that further development of the model requires setting the standard problem of equilibrium refinement. 



An equilibrium in the game can also be  characterized  by   equilibrium partitions. 
\begin{definition}[\textbf{equilibrium coalition structures or partitions}]
A set of  partitions $\{ P^*(K)\}$, $\{ P^*(K)\} \subset \mathcal{P}(K)$, of a game $\Gamma(K)$, is a set of equilibrium partitions,  if it is induced by an equilibrium strategy profile $\sigma^*(K)=\Big(\sigma^*_i(K)\Big)_{ i \in N}$. 
\end{definition}
Robustness of an equilibrium to an increase in $K$ is addressed in Section 8. 


\section{Discussion}



At the moment there are two competing game theories: the non-cooperative and the cooperative. In early 60-s Aumann (1960) wrote that  
`the non-cooperative game differs from the cooperative   game chiefly in that the use of correlated strategy vectors that are not also mixed strategies vectors are forbidden'. 
The current   paper demands another description of the difference: 
the non-cooperative  game theory is based on   mapping  of a strategy profile  of all players  (a subset from $\mathbb{R}^N$) into a  payoff profile for all players  (a  bounded subset from $\mathbb{R}^N$). The cooperative game theory is based on   mapping\footnote{ a characteristic function} of  subsets of integer numbers (a subset from  $\mathbb{N}$) into a subset of real numbers (a   subset from $\mathbb{R}$).   

The    model of this paper explicitly assumes  that every  self-interested agent   produces a relevant impact  for all other players, independently from their coalition location.  The resulting impact can be  either positive or negative. This idea requires  modification  of the above to be explained missed in  existing non-cooperative  and   existing cooperative game theories.  

A possible  irrelevance of  cross-coalition externalities   leaves a  space for a co-existence of  the  non-cooperative and  the cooperative game theories.

Inedequacy  of   cooperative game theory to study  coalitions and coalition structures was  noted  earlier by many authors. Maskin (2011) wrote that `features of cooperative theory are problematic because most applications of game theory to economics involve settings in which externalities are important, Pareto inefficiency arises, and the grand coalition does not form'.\footnote{or 
  that 
`characteristic function ... rules out externalities, situations when a coalition payoff depends on what other coalitions are doing', Maskin, (2016)}.   
 Myerson (p.370, 1991) noted that `we need some model of cooperative behavior that does not abandon the individual decision-theoretic foundations of game theory'. Therefore, there  is a   demand for  a  non-cooperative game to study non-cooperative formation of coalition structures.    
  
 \subsection{Referring to existing literature}
 There are volumes of  literature on   coalitions, and this list of highly respected list   authors is far from complete:  Aumann, Hart, Holt, Kurz, Maschler,  Maskin, Moulen, Myerson, Peleg,  Roth, Serrano, Shapley, Schmeidler, Weber, Winter, Wooders and many others. All  the authors address highly significant specific issues.  However addressing a particular issue   is not the same as a solution for  a general problem,   formation of coalition structures from self-interest  fundamentals.  
 The solution  for a fundamental problem is not  abundance in  literature and is not abundant in current literature. For the purpose of this paper one needs to seek   answers for   certain questions:  
\begin{description}
\item[Problem Identification:] What  are  the specific properties of   a general problem? 
\item[Solution Existence:] Do  existing solutions satisfy  these properties? 
\item[Choice of a tool:] Do existing tools    challenge them enough? \end{description}

\subsubsection{Identification of a problem} There are different views on complexities for non-cooperative formation of coalition structures. There are two opinions    above, two   recent ones are below. 

 Serrano (2014) wrote:  `the axiomatic route find difficulties identifying solutions', and 
 that for studying coalition formation   {`it  maybe worth  to use strategic-form   games, as proposed in the Nash Program'}. 
 Ray and  Vohra (2015) wrote  on complexity and contradictions in existing approaches. They offer a systematic view on the field,  based on `collection of coalitions', or a modified cooperative game theory:
 
`{Yet as one surveys the landscape of this area of research, the  first feature that attracts attention is ``{the fragmented nature of the literature}'' ...
The literature on coalition formation embodies two classical approaches that essentially form two parts of this chapter:
(i) The blocking approach, in which we require the immunity of a coalition arrangement to  blocking, perhaps subcoalitions or by other groups which intersect the coalition in question...
(ii) Noncooperative bargaining, in which individuals make proposals to form a coalition, which can be accepted or rejected...}

{After all, the basic methodologies differ apparently at ``{an irreconcilable level}'' over cooperative and noncooperative game approaches...
 }'

Every  presented view describes only a special part of the general problem    and suggests a partial  solution for it.
Existing diversity of views, `irreconcilable' approaches and `{the fragmented nature}',    enable us to conclude:  the  long lasting  problem of structuring people in multiple groups based on their individual preferences and individual actions is  still not well-identified: a `difficult to identify problem' cannot be easily solved in a consistent way 

 This paper considers the general problem as: how to construct coalition structures with multiple coalitions, intra/inter coalition externalities only from non-cooperative actions of players.  
 
\subsubsection{Existence of a solution}

`A problem cannot be solved if it{'}s bounds are unknown' (A. Tchekmarev\footnote{Personal communication   on engineering design.}).    
The current paper dares to suggest a general identification for the  problem and    a consistent,  natural solution for it.  
  This  game generalizes the concept of  non-cooperative games introduced by   Nash to study formation of coalition structure with  many possible coalitions.

\subsubsection{A tool}
How to deal with the problem,   which has  not been solved and has   `{the fragmented nature}'?  The answer  comes  from Albert Einstein: `The significant problems we have cannot be solved at the same level of thinking with which we created them'.  
  The  current   paper  dares to suggest a new tool, presented in Section 3. 
 
 \subsection{Comparison with existing approaches}
\subsubsection{ A threat} A  threat   as a  tool   for  coalition   formation analysis
 was suggested by Nash (1953), and earlier by von Neumann and Morgenstern.    
   It  starts from  a  strategy profile, a threat, from  \textit{a subset} of   players, possibly allocated in different coalitions. 
   Let this profile  be  a threat to someone beyond this subset.   
The threatening players   may  produce externalities   for each other (and negative  externalities are not excluded!). Then how  can we   describe the externalities for members of the subset if the threat    is  an elementary tool? 
At the same time  there maybe some other player(s) beyond the subset, who may benefit from the threat.  But    they may have no  motivation to join  those who are threatening. For example,  some additional adverse externalities may appear. Such issues should be   included into consideration of coalition structure formation.

There is a parallel argument against using   the   threat  as an elementary tool. Assume there are several agents, who individually cannot make a threat to some  others,   and these small agents are  allocated in different coalitions. A credible threat may come only from many  small players.  Does this mean that the small agents should join together?   They may have their own contradictions to join in one coalition. 
Maybe a formal example will be more illustrative here, but the volume of the paper does is excessive. 

\subsubsection{Usage of the non-cooperative approach}
 The justification for  usage of non-cooperative games as a tool for coalition structure formation, comes from   Maskin (2011) and  the recent   remark from 
 Serrano (2014):   
`it  maybe worth  to use strategic-form   games, as proposed in the Nash Program'.  

There is a difference in the research agenda   of this paper  from  the  Nash Program (Serrano 2004). Nash  program aims to    study a   non-cooperative formation of  \textit{one} coalition, but this paper aims to study non-cooperative formation of     coalition structures, possibly with  more   than one coalition. 
 The constructed finite non-cooperative game allows to  study  what can be a cooperative behavior, when  the individuals
 `rationally expand their individual interests' ( Olson, 1971), see the next section. 
 
 \subsubsection{Novelties of the paper}
 Nash (1950, 1951) suggested to construct a non-cooperative game as a mapping from  a set of strategies into a set  of payoffs,   $$(U_i)_{i \in N} \colon \times_{i \in N} S_i  \rightarrow \mathbb{R}^N,\ \ \  (U_i)_{i \in N} < \infty. $$ 

This paper has two additions in comparison to his paper:   construction of a non-cooperative game with an  embedded coalition structure formation mechanism  and  a parametrization of   games by a number of deviators $K$, $K = 1,\ldots, \#N$:  
 $$  \Big\{\Big(\mathcal{U}_i (S(P))\Big)_{i \in N} \Big\}_{ P \in \mathcal{P}(K)} \colon \times_{i \in N} S_i(K)\buildrel \mathcal{R}(K) \over \rightarrow \{ S(P) \colon P \in \mathcal{P}(K) \} \rightarrow  (\mathbb{R}^N)^{p(N,K)},$$
 where $M$ is a number of coalition structures  from $N$ players with the  restriction $K$ on coalition sizes. 
 The game suggested by Nash becomes a partial case for these games, where coalition structures do not matter and   only one player deviates, $K=1$. 

Every result of  the suggested  game consists of two parts: 1/ a payoff profile for all  players; 2/ an allocation of all players over  coalitions.  Equilibrium in mixed strategies always exists and  may not be efficient like in  traditional non-cooperative games. 

We can take $\mathcal{R}(K)$, $K = 1,\ldots, N$, as a trustworthy social institute.  Application of $\mathcal{R}(K)$ evokes two issues: how to construct an equilibrium (modelled in this paper), and how to achieve  social efficiency equilibrium (beyond the scope of this paper). 
Studying efficiency requires an explicit description for all possible  coalition structure formation rules. 
Then one can define  social desirability in terms of final individual payoffs and allocations of players. This is  the big research program left for the future.

Further development of this research leads to a non-cooperative welfare theory with reconsideration of the first and   second welfare theorems and Arrow{'}s impossibility theorem. These projects can be applied in sociology, management, social and political theories, where disagreement on rules of coalition structure formation can be interpreted as a social or a political conflict. All such issues are left for the future.
 The author believes that   this  new  research direction  will  
   enrich non-cooperative  game  theory and provide non-cooperative fundamentals   for social and mechanism design.  

 
 \subsubsection{Differences from cooperative game theory concepts}
 There are volumes of great literature on the cooperative game theory, with the  famous equilibrium concepts:   the strong  Nash equilibrium (sometimes it is considered as  the non-cooperative equilibrium concept),  the Shapley value, the coalition-proof equilibrium,  the  nucleolus,  the   kernel, the  bargaining set and some other.
 
They  share a list of   common differences with the   equilibrium concept introduced here:    every player does not make an individual choice and does not obtain an individual payoff. In equilibrium there are  no intra- and inter- coalition externalities for every player.   A coalition value, (in terminology of   cooperative game theory,   a sum of individual payoffs of  all coalition members), does not  depend on a whole coalition structure. 
 
There are other more specific differences with every existing cooperative game theory  equilibrium concept.    Differences from the core approach of Aumann (1960) are clear:  no restriction that  only one group deviates, no restriction on the direction of a deviation (inside or outside), and a  construction of individual payoffs from a strategy profile of all players. There is no need to assume transferable/non-transferable utilities for players.
  
 The approach does not need to use the blocking coalition approach, which   does not  study  simultaneous deviation of more than one coalition, and ignores  externalities between  deviators in the coalition and the original set of players. Construction of a sequential  game with coalition structure formation  will follow in the next paper.\footnote{M.Kelbert, O.Pushkareva, D. Levando  Formation of coalition structures as an extended form game, mimeo.}
 


The  role of a central planner in the current paper is different from the one introduced by Nash, who `argued that cooperative actions are the result of some process of bargaining' Myerson (p.370, 1991).  
The central planner offers  a predefined family  of coalition structure formation mechanisms, a family of eligible partitions and a family of rules  to construct these partitions from individual strategies of players.

Based on the properties  of the game we can  propose a criterion of stability  of coalition structures and  study self-enforcing properties of an equilibrium. These properties cannot be designed within existing game theories.

 \section{Formal definition of cooperation}

In   the examples from  Section 2  we  have seen that  efficiency does not imply cooperation. 
 This section formalizes cooperation as an allocation of players in one coalition independently from Pareto efficiency.   A suggested cooperation criterion is  sufficient and can be relaxed in many ways.  



\begin{definition}[complete cooperation in a coalition] 
There is  a non-cooperative game $\Gamma(K)$ with an equilibrium $\sigma^*(K)=(\sigma^*_i)_{i \in N}$.  A  set of players $g$,  ``{cooperate completely    in   the coalition   $g$}''   if for every player $i \in g$ there are:
\begin{description}
\item[\textit{ex ante:}] for every $i$ in $g$,  a  coalition $g$ always belongs to an individually chosen coalition structure, $P_i$ , i.e. if $s_i$ is chosen by $i$, and $s_i \in S_i(P_i)$,  then $g \in P_i$. \footnote{However coalition structures chosen by different players maybe different.}   

\item[\textit{ex post:}] every  realized equilibrium coalition structure     contains $g$, i.e.  $g \in  \forall P^*$, where $P^*$ is a formed equilibrium partition of $\Gamma(K)$,
\end{description}
\end{definition}

Cooperation is defined for a game $\Gamma(K)$ with a fixed $K$.  If $K$ increases     the cooperation may evaporate. 
 
 
 After the game is over the coalition $g$ always belongs to every  final equilibrium coalition structure, \textit{disregarding allocation of players in other coalitions}.  A final partition may  differ from an individual  choice, but  will contain the desirable coalition. 

\section{Application: Bayesian games}
Equilibrium mixed strategies may exist inside one coalition, inducing intra-coalition externalities. To show   this, we reconsider  the standard Battle of the Sexes (BoS) game. 

There are  two players, Ann and Bob. Each  has  two options:  to be together   or to be alone. In every option each can choose   where to go:  Boxing or Opera. Every player has four strategies, in total 16 outcomes. Every outcome consists of a payoff profile and a partition (or a coalition structure). Assume both players have preferences over coalition structures: they prefer to be together, as in the game above with  two extrovert players. 

The rules of coalition structure formation mechanism are:
\begin{enumerate}
\item If they both choose to be together, i.e. both choose  the coalition structure $P_{together}=\{Ann, Bob \}$ then:
\begin{enumerate}
\item if  they choose the same  action  (i.e. both chooses Boxing or both chooses Opera), then they they both made the same choice;
\item otherwise they  do not go anywhere, but enjoy just being together;
\end{enumerate}
\item if at least one of them chooses to remain alone, i.e. chooses a partition   $P_{alone}=\{\{Ann\},\{ Bob \}\}$, then each goes alone to where she/he chooses, maybe to   a different Opera or to   a different Boxing match.
\end{enumerate}

  {\tiny{\begin{table}[htp]
\caption{Payoff for  the  Battle of the Sexes game.  {B is for Boxing, O is for Opera. If the players choose to be together, and it is realized, then   each obtains   an additional fixed payoff $\epsilon >0$ in any outcome.}
} 
\begin{center}
\begin{tabular}{|c|c|c||c|c|} \hline
 & $B_{Bob,alone}$ & $O_{Bob,alone} $& $B_{Bob,together} $& $O_{Bob,together} $ \\  \hline
$B_{Ann,alone}$ &\begin{tabular}{c}   $(2;1)^{*}$  \\  $ \{\{1\},\{ 2\} \}$  \end{tabular} & \begin{tabular}{c} (0;0)  \\ $  \{\{1\},\{ 2\} \}$  \end{tabular} &\begin{tabular}{c}  (2;1) \\  $  \{\{1\},\{ 2\} \}$   \end{tabular} &  \begin{tabular}{c} (0;0)\\  $  \{\{1\},\{ 2\} \}$   \end{tabular} \\ \hline
$O_{Ann,alone}$ &  \begin{tabular}{c} (0;0) \\$  \{\{1\},\{ 2\} \}$ \end{tabular} & \begin{tabular}{c}  $(1;2)^{*}$ \\ $  \{\{1\},\{ 2\} \}$ \end{tabular}& \begin{tabular}{c}  (0;0) \\ $ \{\{1\},\{ 2\} \}$ \end{tabular} &  \begin{tabular}{c} (1;2) \\   $  \{\{1\},\{ 2\} \}$ \end{tabular} \\ \hline \hline
$B_{Ann,together}$ &  \begin{tabular}{c} (2;1) \\ $  \{\{1\},\{ 2\} \}$ \end{tabular} &  \begin{tabular}{c} (0;0) \\ $ \{\{1\},\{ 2\} \}$ \end{tabular} &   
\begin{tabular}{c} 
$(2+\epsilon ;1+\epsilon )^{**} $
  \\ $ \{1, 2 \}$ 
  \end{tabular} & \begin{tabular}{c}  $ (\epsilon; \epsilon)$  \\  $  \{1,2\}$ \end{tabular} \\ \hline
$O_{Ann,together}$ & \begin{tabular}{c}  (0;0) \\ $ \{\{1\},\{ 2\} \}$ \end{tabular} & \begin{tabular}{c}   (1;2) \\ $ \{\{1\},\{ 2\} \}$ \end{tabular}& 
\begin{tabular}{c}  $(\epsilon;\epsilon) $ \\ $ \{1,2\}$ \end{tabular} & 
\begin{tabular}{c}  $(1+\epsilon;2+\epsilon)^{**}$ \\  $ \{1, 2\}$ \end{tabular}  \\ \hline  
 \end{tabular}
\end{center}
\label{default2}
\end{table}%
 }}
 
 Formally the coalition structure formation rules  are: 
\begin{multline*}
 \mathcal{R}(K=1) \colon S(K=1) \mapsto S(P=\{ \{Ann \} ,\{Bob \} \}), \\  \forall s \in S_i(K=1)=
 \{O_{Ann,alone}, B_{Ann,alone } \} \times \{O_{Bob,alone},B_{Bob,alone }\}
  \end{multline*}
 and
 
 $$
 \mathcal{R}(K=2) \colon S(K=2) \mapsto \begin{cases}
  S(P=\{Ann, Bob\}) = \Big\{ s \colon  \\  s \in  \{O_{Ann,together}, B_{Ann,together} \} \times \{O_{Bob,together},B_{Bob,together }\}  
  \\ S\Big(P=\{\{Ann\},\{Bob \}\}\Big) = \\ S(K=2) \setminus S(P=\{Ann,Bob \} \})  \ \ 
 \mbox{otherwise}
 \end{cases}.
 $$

 Table  \ref{default2} corresponds to the game with $K=2$, where the game for $K=1$ is a nested component. 
  If Ann and Bob play the game with $K=1$,   a single coalition structure $\{\{ Ann\} , \{ Bob\} \}$ game, then the payoffs for this game are  in the two-by-two top-left corner of Table  \ref{default2}, and this is the standard game.   
 If Ann and Bob are together, then each obtains an additional payoff $\epsilon$,  and the corresponding cells make the two-by-two  bottom-right corner of Table \ref{default2}.
 
 BoS game with one deviator ($K=1$)  has two equilibria in pure strategies and  one in mixed. For  $K=2$   only mixed strategy equilibrium survives, as  two player can deviate simultaneously. 
 Mixed  equilibrium for  $K=2$  differs from one for $K=1$: another domain of pure strategies makes  another  coalition structure and  another  payoff profile. Mixed strategies equilibrium  for Ann is: $\sigma^*(B_{Ann,together})=(1+\epsilon)/(3+2\epsilon)$, 
 $\sigma^*(O_{Ann,together})=(2+\epsilon)/(3+2\epsilon)$.

Equilibrium mixed strategies may appear   for every $K=1,2$.  But for $K=2$ equilibrium coalition structure is  the grand coalition, where individual payoffs fluctuate. More of that, for $K=2$ a value of  the grand coalition holds constant for equilibrium  mixed strategies. Such games  are not described in available literature.

 \section{Application: Stochastic games}
 Shapley  (1953) defined  stochastic games as `the play proceeds by steps from position to position, according to transition probabilities controlled jointly by the two player'. 
This section demonstrates how  to construct this type of games. 

The  example below  differs from   the one above in two respects. A set of  equilibrium mixed strategies induces more than one equilibrium coalition structure.   
Another difference is that   an individual payoff of a coalition  depends on a coalition structure to which it belongs. This  game is the closest to studying group/paintball team examples in the introduction of the paper.

\subsection{Corporate lunch game}

There is a set of four identical players $N=\{A,B,C,D\}$.  
A coalition structure is an allocation of   players over no more than four separate tables. \footnote{Possibly empty tables do not matter, and tables  can not be moved.} 
A strategy of $i$   is a coalition structure, or an allocation of all players across tables during  a lunch.    Different from previous examples there are no strategies inside coalitions, a pure strategy is a coalition structure.

A rule of coalition structure formation is that any coalition in any coalition structure can be formed only from a  unanimous agreement from members of this coalition. 

A player has preferences over coalition structures:  she/he prefers to eat with another one, but also  prefers that all  other players eat individually. Possible motivation could be a   dissipation of  information or gossips.  From the other side, if one eats alone he/she is hurt by     formed coalition of others.

 If   individual choices  of coalition structures are inconsistent, then  we  apply the rule of coalition structure formation.  
The first column of Table  \ref{lunch} contains formed  coalition structures. 

\begin{table}[htp]
\caption{Office lunch game: strategies and payoff profiles. Full set of equilibrium mixed strategies are indicated  only for player $A$.}
 \begin{center}
\begin{tabular}{|c|c|c|c|}
\hline num&
Coalition structure &  Payoff profile $U_A,U_B,U_C,U_D$ & \begin{tabular}{c} \small{Coalition values } \\ \small{as in }  \\  \small{ cooperative game theory} \end{tabular}  \\
\hline
$1^*$& $\{ A,B\}, \{ C\}, \{D\}$:  $ \sigma^*_A= \sigma^*_B=1/3$ & (10,10,3,3)  & $20_{A,B}$, $3_{C}$, $3_{D}$ \\
$2^*$& $\{ A,C\}, \{ B\}, \{D\}$:  $ \sigma^*_A=\sigma^*_C=1/3$ & (10,3,10,3)  & $20_{A,C}$, $3_{B}$, $3_{D}$ \\
$3^*$& $\{ A,D\}, \{ C\}, \{B\}$:  $ \sigma^*_A=\sigma^*_D=1/3$ & (10,3,3,10)  & $20_{A,C}$, $3_{C}$, $3_{B}$ \\
4 & $\{ A\}, \{B\}, \{ C, D\}$ & (3,3,10,10)  & $3_{A}$, $3_{B}$, $20_{C,D}$\\
5 & $\{ A\}, \{D\}, \{ B,C\}$ & (3,10,10, 3)  & $3_{A}$, $3_{D}$, $20_{C,B}$\\
6 & $\{ A\}, \{C\}, \{ B, D\}$ & (3,10,3,10)  & $3_{A}$, $3_{C}$, $20_{B,D}$\\
7& $\{ A\} ,\{B\}, \{ C\}, \{D\}$ & (3,3,3,3) & $3_{A}$, $3_{B}$, $3_{C}$, $3_{D}$\\
8& $\{ A, B\}, \{ C,D\}$ & (3,3,3,3) & $6_{A,B}$, $6_{C,D}$\\
9& $\{ A,C\}, \{ B,D\}$ & (3,3,3,3) & $6_{A,C}$, $6_{B,D}$\\
10 & $\{ A, D\}, \{ B,C\}$ & (3,3,3,3) & $6_{A,D}$, $6_{B,C}$\\
11  & \mbox{all other with $K=3,4$} & (0,0,0,0) &  $ = 0$\\
\hline
\end{tabular}
\end{center}
\label{lunch}
\end{table}%
 
 Payoffs in   Table \ref{lunch} are organized to satisfy  the individual preferences described above.  Take column 4 and compare payoff for the coalition $\{A,B \}$ in different coalition structures: line 1, $\{\{A,B \}, \{ C\} , \{ D\}\}$ and line 8, $\{\{A,B \}, \{ C, D\}\}$. 
 Player $C$ prefers to eat with $D$, but also prefers  $A$ and $B$ do not eat together. This is described as that    
 payoffs of the coalition $\{ C,D\}$ depend on the whole coalition structure. \textit{This result is impossible in existing cooperative game theory:  a coalition value depends on allocation of all players.}
 
  
  Let   maximum number of players in a coalition be  2, $K=2$.  In equilibrium every player  chooses  those coalition structures, where she/he is with someone, while others  are separate.  Equilibrium mixed strategies  are indicated only for player $A$. 
  
Several different coalition structures can form in an equilibrium. Each of them can be considered as a state of a stochastic game.

An increase in $K$,    maximum coalition size,   decreases payoffs for all players. An increment of   $K$ from $2$ to $3$ and then to $4$ does not change an equilibrium in mixed strategies. Thus the equilibrium for $K=2$ is robust to greater $K$s.  This issue is addressed Section 8.

 
 
 \subsection{A formal definition of a stochastic game of coalition structure formation}
 
 Let $\Gamma(K)$ be a non-cooperative game as defined above.
 
 \begin{definition}
 A game $\Gamma(K)$ is a stochastic game if a set of equilibrium partitions is bigger than one, $\# \Large( \{P^*\}(K)\Large) \ge 2$, where a state  is an equilibrium partition $P^*$, $P^* \in \{P^* \}(K)$.
 \end{definition}
 
  Studying properties of stochastic games with non-cooperative coalition structure formation   is left for the  future.
  

 \section{Application: non-cooperative criterion for stability}
The  Nash equilibrium, as defined in Nash (1950, 1951), does not specify how a player may deviate within a group. This leads to  cases, when  desirable outcomes cannot be supported by    Nash equilibrium, while an outcome intuitively  requires   studying possible  deviations of more than one player. Sometimes such points are called `focal points' (Schelling).
  
  Aumann (1990)   demonstrates absence of a  self-enforcement property of Nash equilibrium for  a focal point in `a stag and   hare game'.  The key problem is that standard Nash equilibrium does not have a tool to study deviations of more than one player.
 An  example below explains how to change the  basic game to make joint hunting for a stag  become an equilibrium focal point robust to a deviation of two players.  Then  the paper presents a  non-cooperative   criterion  to measure stability of a coalition structure. 
   
 There are two hunters $i=1,2$.
 If  they can hunt (only) separately, then $K=1$, and  the only available partition is $P_{alone}=\{ \{ 1\}, \{ 2\}\}$.   An individual strategy  set  of $i$ is
$$S_i(K=1)=\Big\{({P_{alone}, hare}),({P_{alone}, stag}) \Big\} $$ with a general element $s_i$. Every $s_i$ consists of two terms: who is the hunting partner (oneself) and what is the  animal to hunt. For example, a strategy $s_i=({P_{separ}, hare})$ is interpreted as  a strategy for  $i$   to  hunt alone for a hare. A set of   corresponding strategies for the game with $K=1$ is $$S(K=1)=S_1(K=1) \times S_1(K=1) .$$  

For a game with $K=2$ every hunter  can additionally choose  to hunt in a company of two and as before to choose   a target for hunting: a hare or a stag. 
 A set  of strategies of $i$ is
 $$S_i(K=2)=\Big\{({P_{alone}, hare}),({P_{alone}, stag}), ({P_{together}, hare}),({P_{together}, stag}) \Big\}. 
 $$  
 A set of strategies of the game is a direct (Cartesian) product, $$S(K=2)=S_1(K=2) \times S_2(K=2).$$

  We do not rewrite the rules for coalition structure formation, as they are the same as  in the Battle of the Sexes  game above. The difference is   in  renaming   strategies. 
  Every player    knows, which game is played, either with $K=1$ or with $K=2$,  and this is the common knowledge.

Payoffs for the   games $\Gamma(K=1)$ and $\Gamma(K=2)$ are   in Table \ref{harestag}. Some  payoff   profiles  have a special interpretation: $(8;8)$ -  every hunter obtains a hare, $(4;4)$ - two hunters obtain a hare for two, $(100;100)$ -  both hunters obtain one stag. The point  $(100;100)$ is a desirable outcome and a focal point, which  cannot be supported as the Nash equilibrium with one deviator Aumann (1990).

For a  game with $K=1$    maximum achievable payoff is $(8,8)$, when each hunts individually for a hare. An   equilibrium strategy profile  is
$$s^*(K=1)=\Big((P_{alone}, hare),(P_{alone}, hare) \Big)$$ with the payoff  profile $(8;8)$. 
 In the game $\Gamma(K=1)$   the  players  cannot reach the  efficient outcome $(100,100)$. It is available only if   both hunters can decide to hunt together.
  This focal point (in terminology of Schelling) can  be  reached only within the game $\Gamma(K=2)$, but not for the game $\Gamma(K=1)$. This is the explanation for the problem posed  by Aumann:  there can be an attractive outcome of a game, but  it   cannot be described  in terms of a Nash equilibrium of a traditional non-cooperative game. Deviation of two players demands additional tools  for the model. 
  
  An   equilibrium strategy profile is one  when  hunters can operate together: 
$$s^*(K=2)=\Big((P_{together}, hare),(P_{together}, hare) \Big)$$ and reach the payoff  profile $(100;100)$. 
The focal point  $(100;100)$ is an equilibrium for the game $\Gamma(K=2)$, but it is unfeasible for  the game $\Gamma(K=1)$.


    \begin{table}[htp]
\caption{Expanded stag and hare game}
\begin{center}
\begin{tabular}{|l|c|c|c|c|}
\hline 
&\small{\mbox{\negthickspace  $P_{alone}$, hare}} &\small{\mbox{\negthickspace  $P_{alone}$, stag}} & \small{\negthickspace\mbox{ $P_{together}$,hare}} & \small{\negthickspace \mbox{ $P_{together}$, stag}} \\
\hline
\small{\mbox{$P_{alone}$,hare}} &  $(8;8)^{*}$; $\{ \{ 1\},\{ 2\}\}$ & (8;0); $\{ \{ 1\},\{ 2\}\}$  & (8;8);$\{ \{ 1\},\{ 2\}\}$ & (8;0); $\{ \{ 1\},\{ 2\}\}$ \\
\small{\mbox{$P_{alone}$, stag} }  & (0;8); $\{ \{ 1\},\{ 2\}\}$ & (0;0); $\{ \{ 1\},\{ 2\}\}$& (0;8); $\{ \{ 1\},\{ 2\}\}$& (0;0); $\{ \{ 1\},\{ 2\}\}$\\
\hline
\small{\mbox{$P_{together}$, hare} } & (8;8); $\{ \{ 1\},\{ 2\}\}$ & (8;0); $\{ \{ 1\},\{ 2\}\}$  & (4;4); $\{ 1, 2\}$ & (8;0); $\{ 1, 2\}$\\
 \small{\mbox{$P_{together}$, stag}}   & (0;8); $\{ \{ 1\},\{ 2\}\}$& (0;0) ;$\{ \{ 1\},\{ 2\}\}$ & (0;8);$\{ 1, 2\}$ & $(100;100)^{**}$; $\{ 1, 2\}$\\
 \hline
\end{tabular}
\end{center}
\label{harestag}
\end{table}%

 If it is  uncertain, which game is played,  either $\Gamma(K=1)$ or $\Gamma(K=2)$, then players will randomize between two strategies:
 $(P_{alone}, hare)$ and $(P_{separ}, stag)$ , and the game become   a stochastic game.
 
 \subsection{Criterion of  coalition structure (a partition) stability}
There is a    nested family  of games $$\mathcal{G} = \{\Gamma(K=1), \ldots, \Gamma(K), \ldots \Gamma(K=N) \} \colon$$ such that $$\Gamma(K=1) \subset \ldots \subset  \Gamma(K) \subset \ldots \subset  \Gamma(K=N).$$ 
A list of corresponding   mixed strategies equilibria is  $$\Big( \sigma^*(1), \ldots, \sigma^*(K), \ldots, \sigma^*(K=N)\Big),$$ where  $ \sigma^*(K)=( \sigma^*_i(K))_{i \in N}$. Resulting  equilibrium  coalition structures (partitions) are :
 $$\Big(\{ P^*\}(1), \ldots, \{ P^*\}(K),\ldots, \{ P^*\}(K=N) \Big),$$ where for every $K=1,\#N$ there is $\{ P^*\}(K)  \subset \mathcal{P}(K)$.

The family of games has an equilibrium expected  payoff profiles:

$$
\Big( (EU^{\Gamma(1)}_i)^*_{i \in N}, \ldots,( EU^{\Gamma(K)}_i)^*_{i \in N}, \ldots,( EU^{\Gamma(K=N)}_i)^*_{i \in N} \Big) 
, $$ where $( EU^{\Gamma(K)}_i)^*_{i \in N} \equiv ( EU^{\Gamma(K)}_i (\sigma^*))_{i \in N}$. 

Take a game $\Gamma(K_0)$ from $  \mathcal{G}$ with   maximum coalition size $K_0$. It  has $\sigma^*(K_0)$ as an equilibrium mixed strategy profile. The question is: what is the condition when   an equilibrium strategy profile  does not  change with an increase in   maximum  coalition  size $K$?   

  The suggested  criterion is based on the idea that a set of mixed strategies should not change with an increase in  $K$, $\sigma^*(K_0)=\sigma^*(K_0+1) = \ldots =  \sigma^*(N)$, with the restriction
  $Dom (\sigma^*(K_0))= Dom (\sigma^*(K_0+1)) = \ldots =  Dom(\sigma^*(N))$, where $Dom$ is a domain of   equilibrium mixed strategies. The same in other words,   if a game $\Gamma(K_0)$ is  substituted  sequentially by   games $\Gamma(K+1)$ upto $\Gamma(N)$, then an equilibrium strategy profile will   not change in domain and in range.
   The next criterion   is a sufficient criterion for coalition structure stability. 
 \begin{definition}
\textbf{Coalition structure (partition) stability criterion}   for a game $\Gamma(K_0)$ is a  maximum coalition size $K^*$, when an equilibrium  still holds true, i.e. for all $i \in N$ there is a   number $K^*$ such that
   \begin{enumerate}
   \item $$K^* =   \max_{K = K_0,\ldots,  N \atop{\Gamma(K_0) \ldots, \Gamma(K=N)}}  \Big\{ EU^{\Gamma(K_{0})}_i\Big(\sigma^*_i (K_{0}),
  \sigma^*_{-i}(K_{0})\Big) \ge \\  EU^{\Gamma(K)}_i\Big(\sigma^*_i (K),
  \sigma^*_{-i}(K)\Big)  \Big\},$$
  \item 
 $Dom \mbox{ }  \sigma^*(K^*)  = Dom \mbox{ }  \sigma^*(K_0)$
 \end{enumerate}
where $\sigma^*(K_0)$ is an equilibrium in \underline{the} game $\Gamma(K_0)$,  while $\sigma^*(K)$ is an equilibrium in \underline{a} game $\Gamma(K)$, $K = K_0,\ldots,  N$, and  $Dom$  is a  domain of equilibrium mixed strategies set.

\end{definition} 

The definition is operational, it can be constructed directly from a definition of a family of games $\mathcal{G}$.
It guarantees stability of both payoffs and partitions. In Prisoner{'}s Dilemma example (where both players are ``non-extroverts'', ``non-introverts'') an increment in $K$    increases a number of equilibrium coalition structures without a change in payoffs. This increases a number of equilibria, without  rejecting one for $K=1$. The criterion suggests $K^*=2$  for this game. 
The lunch game  (Section 7) example has robust equilibrium for $K^*=2$. 
Some applications may require weaker forms of the  criterion. 


 In the extended version of stag and hare game   an increase in $K$ changed the equilibrium. 
The same took place in a variation of Battle of the Sexes game. 
However this did not happen in  the Corporate Lunch game.

The proposed criterion may serve as a measure of trust to an equilibrium or as a test for self-enforcement   of an equilibrium.
This criterion can be  applied to study opportunistic behavior  in coalition partitions. If  players  in a coalition $g$ of a game $\Gamma(K_1)$  have perfect cooperation (Section 5), this does not mean that in a wider game  $\Gamma(K_2)$,  $K_1 < K_2$,  they will still cooperate. 

Studying stability is tightly connected to studying focal points.  The suggestion of the criterion is only the first step to  the big field of   research. 
\section{Conclusion}

 Traditionally social sciences are interested in structuring people in multiple groups based on their individual preferences. The   paper presents a family of non-cooperative finite games for coalition structure formation motivated by self-interest benefits. Every game in a family has an equilibrium in mixed strategies. Examples demonstrate how   the approach can be applied to study some popular games.  The paper  introduces a non-cooperative criterion  to measure stability of coalition structures   based on   maximum number of deviators. 
 
 The novelty of the paper is not in mathematical, but in the  structural result. The suggested family of games can help to study  new issues, earlier beyond the scope of research, for example, focal points.
Development of the model for  repeated and network games  are in   other forthcoming papers from the author and his colleagues.

\end{document}